# Novel General Active Reliability Redundancy Allocation Problems and Algorithm

Wei-Chang Yeh
Integration and Collaboration Laboratory
Department of Industrial Engineering and Engineering Management
National Tsing Hua University
yeh@ieee.org

**Abstract:** The traditional (active) reliability redundancy allocation problem (RRAP) is used to maximize system reliability by determining the redundancy and reliability variables in each subsystem to satisfy the volume, cost, and weight constraints. The RRAP structure is very simple, that is, redundant components are parallel in each subsystem, and all subsystems are either connected in series or in a bridge network. Owing to its important and practical applications, a novel RRAP, called the general RRAP (GRRAP), is proposed to extend the series-parallel structure or bridge network to a more general network structure. To solve the proposed novel GRRAP, a new algorithm, called the BAT-SSOA3, used the simplified swarm optimization (SSO) to update solutions, the small-sampling tri-objective orthogonal array (SS3OA) to tune the parameters in the proposed algorithm, the binary-addition-tree algorithm (BAT) to calculate the fitness (i.e., reliability) of each solution, and the penalty function to force infeasible back to the feasible region. To validate the proposed algorithm, the BAT-SSOA3 is compared with state-of-the-art algorithms, such as, particle swarm optimization (PSO) and SSO, via designed experiments and computational studies.

## 1. INTRODUCTION

To date, two main strategies have been implemented to improve system/network reliability: 1) increase each component reliability and/or 2) adapt redundant components in subsystems [1]. The traditional (active) reliability redundancy allocation problem (RRAP) using redundant components in parallel is a well-known optimization problem for enhancing the system reliability effectively and economically [2-10]. Hence, RRAP is common in real-life applications, such as the Internet of Things [11] and wireless sensor networks [12]. RRAP determines the redundancy variables (integer variables) and reliability variables (floating variables) in each subsystem to optimize the system reliability in satisfying the resource consumption nonlinear constraint, such as the volume, cost, and weight [1-12].

The main structure in the traditional RRAP is series-parallel, that is, subsystems are in series and the components are parallel to each subsystem [1-12]. In fact, series-parallel networks are not common, and instead, it is more reasonable that most network structures depend on their practical applications [1]. To the best of our knowledge, no research has discussed the general RRAP (GRRAP) by considering the general subsystem structure.

After generalizing the network structure from series-parallel to any structure, the major challenge is calculating the network reliability [1], which is NP-hard [13]. System/network reliability is the success probability that the entire system/network is still functioning. In traditional RRAP, the reliability of each network is simple and can be easily calculated manually [1-12]. However, it turns to a binary-state network reliability problem after generalizing to any type of network structure. Hence, to calculate the binary-state network reliability, the binary-addition-tree algorithm (BAT) proposed by Yeh is adapted here to achieve the abovementioned goal [14].

After applying the algorithm to calculate the binary-state reliability, the next step is to determine the number and reliability of the redundant components. The variable representing the number of redundant components is called redundancy variable, and that representing the reliability of redundant components is called reliability variable. The original RRAP is a mixed-integer nonlinear programming problem in which each redundancy variable is an integer, each reliability variable is a floating-point number, and the objective and constraints are nonlinear [7, 15].

The traditional RRAP is NP-complete [15], and scholars have considered numerous algorithms to solve RRAP from various viewpoints [1-12]. These algorithms are mainly categorized as follows.

1. Exact-solution algorithms, such as column generation approach [16, 17] and branch-and-bound [18], are only available for small-sized RRAPs to obtain optimal solutions in an acceptable amount of time.

2. Heuristic algorithms [8], such as surrogate constraints algorithm [19], require derivatives that

are not easily derived and extensive computational efforts to obtain near-optimal solutions for medium-sized RRAPs.

3. Artificial Intelligent (AI) algorithms, such as simplified swarm optimization (SSO) [5, 22, 23], particle swarm optimization (PSO) [20], hybrid swarm optimization (HSO) [2, 3, 23], genetic algorithms (GA) [7], artificial immune algorithm (IA) [4], artificial bee colony algorithm (ABC) [6, 21], and PSSO (combining both PSO and SSO) [23], are used to reduce the computational burden for obtaining approximated solutions for larger RRAPs.

Among these algorithms, AI algorithms have gained much attention as they can improve the efficiency and solution quality of existing algorithms, resulting in good-quality solutions for solving large practical problems [2, 3, 4, 5, 6, 20, 21, 22, 23].

Among various AI algorithms, SSO first proposed by Yeh [20, 23] is attractive because it is easy to understand and implement, efficient, in most cases, effective in obtaining high-quality solutions, simple to code using different computer languages, and sufficiently flexible to be used with other algorithms (hybrid algorithms) [2, 3, 5, 22, 23]. Hence, the goal of this study is to extend the traditional RRAP to GRRAP and solve the novel GRRAP based mainly on BAT and SSO [1].

The remainder of this paper is organized as follows. Section 2 introduces the traditional RRAP, BAT, SSO, and penalty functions. Section 3 presents the proposed BAT-connected vectors, novel small-sampling tri-objective orthogonal array (SS3OA), and the proposed new BAT-SSOA3 hybrid with the boundary condition differentiated from the cost constraint, BAT, and SS3OA. In Section 4, a comprehensive comparative study on the performance of the proposed BAT-SSOA3 and existing methods adapted from the traditional RRAP is presented. Finally, the conclusions are presented in Section 5.

## 2. RRAP, BAT, SSO, AND PENALTY FUNCTION

The traditional RRAP s generalized to GRRAP, and the proposed BAT-SSOA3 is based on BAT, SSO, and penalty functions to solve the GRRAP. Hence, before presenting the proposed BAT-SSOA3, the required notations, the traditional RRAP, and the reasons for generalizing RRAP, basic BAT, simple

SSO, and a common penalty function are introduced briefly in the following subsection.

## 2.1 Notations

The notations used in the study are listed below.

$N_{var}$ ：number of subsystems/redundancy variables/reliability variables.

$N_{sol}$ ：number of solutions.

$N_{gen}$ ：number of generations

$\rho_I$ ：number generated from interval $I$ randomly.

$n_i$ ：redundancy variable in subsystem $i$ for $i = 1, 2,…, N_{var}$.

$r_i$ ：reliability variable in subsystem $i$ for $i = 1, 2,…, N_{var}$.

$\mathbf{n}$ ：$\mathbf{n} = (n_1, n_2, ..., n_{Nvar})$

$\mathbf{r}$ ：$\mathbf{r} = (r_1, r_2, ..., r_{Nvar})$.

$X_k$ ：$X_k = (\mathbf{n}_k, \mathbf{r}_k) = (x_{k,1}, x_{k,2}, …, x_{2k,Nvar})$ for $k = 1, 2, …, N_{sol}$.

$P_k$ ：$pBest$ of the $k^{th}$ solution $P_k = (p_1, p_2,… p_{2·Nvar})$ for $k = 1, 2, …, N_{sol}$.

$G$ ：$gBest$ $G = (g_1, g_2,… g_{2·Nvar})$ which is the best solution among all solutions.

$P_{gBest}$ ：$P_{gBest} = G$

$R_s(\mathbf{n}, \mathbf{r})$ ：system/network reliability under $\mathbf{n}$ and $\mathbf{r}$.

$R_p(\mathbf{n}, \mathbf{r})$ ：penalized system/network reliability under $\mathbf{n}$ and $\mathbf{r}$.

$g_v(\mathbf{n}, \mathbf{r})$ ：volume constraint under $\mathbf{n}$ and $\mathbf{r}$.

$g_c(\mathbf{n}, \mathbf{r})$ ：cost constraint under $\mathbf{n}$ and $\mathbf{r}$.

$g_w(\mathbf{n}, \mathbf{r})$ ：weight constraint under $\mathbf{n}$ and $\mathbf{r}$.

$F(\bullet)$ ：fitness function of solution $\bullet$.

$\bullet_{ub}$ ：upper bound of $\bullet$, e.g., $V_{ub}$, $C_{ub}$, $W_{ub}$

## 2.2 Reliability Redundancy Allocation Problem

RRAP is an important mechanism that gradually emerges in the early stages of planning, designing, and controlling systems/networks [24]. The traditional RRAP has a maximum of five subsystems, that is, $N_{var} \le 5$. $\mathbf{n} = (n_1, n_2, …, n_{Nvar})$ and $\mathbf{r} = (r_1, r_2, …, r_{Nvar})$ are assigned as the vectors formed by redundancy and reliability variables, respectively, where $0 \le r_1, r_2, r_3, r_4, r_5 \le 1$, and $n_1, n_2, n_3, n_4, n_5 = 1, 2, …, 10$ [2, 3, 4, 5]. The general model for RRAP can be defined using mixed-integer nonlinear programming as follows [2, 3, 4, 5, 24]:

Maximize    $R_s(\mathbf{n}, \mathbf{r})$                                    (1)

$$\text{Subject to} \quad g_v(\mathbf{n}, \mathbf{r}) = \sum_{i=1}^{5} w_i v_i^2 n_i^2 \leq V_{ub} \tag{2}$$

$$g_c(\mathbf{n}, \mathbf{r}) = \sum_{i=1}^{5} \alpha_i (-1000 / \ln r_i)^{\beta_i} [n_i + \exp(n_i / 4)] \leq C_{ub} \tag{3}$$

$$g_w(\mathbf{n}, \mathbf{r}) = \sum_{i=1}^{5} w_i n_i \exp(n_i / 4) \leq W_{ub}. \tag{4}$$

Eq. (1) is the objective of maximizing the system reliability $R_s(\mathbf{n}, \mathbf{r})$ of the RRAP. Eqs. (2) – (4) are nonlinear constraints for the volume, cost, and weight, respectively.

The four benchmark RRAPs are depicted in Fig. 1 [2, 3, 4, 5, 24], and the corresponding data, that is, $\alpha_i \cdot 10^5$, $\beta_i$, $w_i v_i^2$, $w_i$, $V_{ub}$, $C_{ub}$, and $W_{ub}$, for the series system shown in Fig. 1(1) and the bridge system shown in Fig. 1(3) are listed in Table 1.

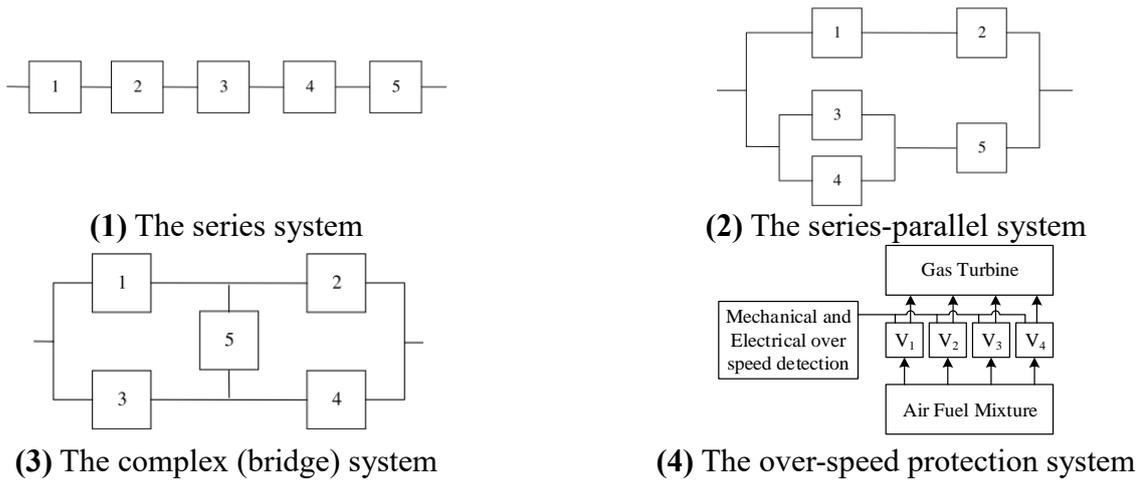

**(1)** The series system

**(2)** The series-parallel system

**(3)** The complex (bridge) system

**(4)** The over-speed protection system

**Figure 1.** Four benchmark RRAPs.

**Table 1.** Data used in Fig. 1(1) and 1(3) [2, 3, 4, 5, 24].

| Subsystem $i$ | $\alpha_i \cdot 10^5$ | $\beta_i$ | $w_i v_i^2$ | $w_i$ | $V_{ub}$ | $C_{ub}$ | $W_{ub}$ |
|---|---|---|---|---|---|---|---|
| 1 | 2.330 | 1.5 | 1 | 7 | 110 | 175 | 200 |
| 2 | 1.450 | 1.5 | 2 | 8 | | | |
| 3 | 0.541 | 1.5 | 3 | 8 | | | |
| 4 | 8.050 | 1.5 | 4 | 6 | | | |
| 5 | 1.950 | 1.5 | 2 | 9 | | | |

The following property was first proposed in [2, 3] to perform a local search to improve the solution quality, including the adjustment of an infeasible solution to a feasible solution.

**Property 1:** Let the values of both $n_i$ and $r_k$ be constants for $i \neq j \in \{1, 2, \dots, N_{var}\}$. $R_s(\mathbf{n}, \mathbf{r})$ is maximized

if [2, 3]

$$r_j = \exp\left( \frac{(-1000) \cdot \sqrt[\beta_j]{\alpha_j \left[ n_j + \exp(\frac{n_j}{4}) \right]}}{\sqrt[\beta_j]{C_{ub} - \sum\limits_{\substack{i=1 \\ i \neq j}}^{Nvar} \alpha_i (\frac{-1000}{\ln r_i})^{\beta_i} \left[ n_i + \exp(\frac{n_i}{4}) \right]}} \right). \tag{5}$$

**Example 1.** Let $C_{UB}$=191, $\mathbf{n}$ = (3, 2, 2, 3, 3), and $\mathbf{r}$ = (0.77946645, 0.87173278, 0.90284951, 0.71148780, 0.78781644). Eq. (17) is implemented to improve the value of the last reliability variable, i.e., $r_5$, in $\mathbf{r}$. After that, we have new $\mathbf{r} = \mathbf{r}^*$ = (0.77946645, 0.87173278, 0.90284951, 0.71148780, 0.7878166527), $R_s(\mathbf{n}, \mathbf{r})$ = .93168229721527107 is increased to $R_s(\mathbf{n}, \mathbf{r}^*)$ = 0.9316823242437, and $g_c(\mathbf{n}, \mathbf{r})$ = 174.999954 is also increased to $g_c(\mathbf{n}, \mathbf{r}^*)$ = 174.99999999999997.

Many comprehensive works on RRAPs can be found in the existing literature, and the details can be found in [2, 3]. As shown in Fig. 1, each benchmark is relatively simple and cannot satisfy many networks for a particular purpose in practical life. Hence, there is a need to extend the traditional RRAP to a general RRAP to meet all types of networks in real life.

The analytical solution for Fig. 1(3) is given by Eq. (6), which shows that it is difficult to calculate the network reliability even for the small-sized benchmark, as shown in Fig. 1(3). Thus, a general RRAP, although important, is difficult to solve [3].

$$R_s(\mathbf{n},\mathbf{r}) = r_1^{n_1} \cdot r_2^{n_2} + r_3^{n_3} \cdot r_4^{n_4} + r_1^{n_1} \cdot r_4^{n_4} \cdot r_5^{n_5} + r_2^{n_2} \cdot r_3^{n_3} \cdot r_5^{n_5} - r_1^{n_1} \cdot r_2^{n_2} \cdot r_3^{n_3} \cdot r_4^{n_4} - r_1^{n_1} \cdot r_2^{n_2} \cdot r_3^{n_3} \cdot r_5^{n_5}$$
$$- r_1^{n_1} \cdot r_2^{n_2} \cdot r_4^{n_4} \cdot r_5^{n_5} - r_1^{n_1} \cdot r_3^{n_3} \cdot r_4^{n_4} \cdot r_5^{n_5} - r_2^{n_2} \cdot r_3^{n_3} \cdot r_4^{n_4} \cdot r_5^{n_5} + 2r_1^{n_1} \cdot r_2^{n_2} \cdot r_3^{n_3} \cdot r_4^{n_4} \cdot r_5^{n_5}. \tag{6}$$

## 2.3 Binary-Addition-Tree Algorithm

Yeh [14] proposed the first BAT algorithm, which was based on binary addition, to find all possible vectors (solutions) in the solution space. $X = (x_1, x_2, \ldots, x_{Nvar})$ was assigned as the vector, and the current coordinate $x_i$ as the state of $a_i$. Note that $x_i$ is either 0 or 1 for $i$ = 1, 2, …, $N_{var}$ in the binary-state network, which is the structure of the general RRAP. The two main principles of implementing the (forward) BAT are as follows:

1. If $x_i = 0$, $x_i$ is changed to 1, and a new vector is obtained. For example, vector $X = (0, 0, 1)$ is updated to $(1, 0, 1)$ if $x_i = x_1 = 0$, that is, 001 is updated to 101 by adding 1 to the first coordinate using binary addition, where $X$ is the current vector and $(1, 0, 1)$ is the new $X$.

2. If $x_i = 1$, $x_i$ is changed to 0, and proceeds to the next coordinate, that is, $x_{(i+1)}$, to repeat the first principles above. For example, $(1, 0, 1)$ is updated to $(0, 1, 1)$, that is, 101 is updated to 011 by adding 1 to the first coordinate using binary addition.

Each vector obtained from the BAT following the first step above corresponds to a subnetwork, and all obtained vectors do not need to be stored, saving computer memory. However, in the proposed BAT-SSOA3, any vector is stored, if there is a path connecting the input and output in its corresponding subnetwork. The BAT pseudocode is given as follows [14, 25, 26]:

**Algorithm: BAT**

**Input:**     A network $G(V, E)$ with the input node 1 and the output node $n$.

**Output:**    All vectors without duplications.

**STEP B0.**  Let $N_{var}$-tuple $X$ be a vector zero and $i = 1$.

**STEP B1.**  If $x_i = 1$, let $x_i = 0$. Otherwise, let $x_i = 1$, $i = 1$, and go to STEP B1.

**STEP B2.**  If $i < N_{var}$, let $i = i + 1$ and go to STEP B1. Otherwise, halt.

STEP B0 initializes $X = 0$ and $i = 1$. STEP B1 implements two principles, and STEP B2 implements the stopping criterion. The time complexity for the BAT is $O(2^{m+1})$ to determine all possible vectors from vector zero to vector one [25], where $m = N_{var}$. Through the use of the pseudocode, the BAT is easy to code, efficient to execute, and effective in saving computer memory.

For example, Fig. 2. is an activity-on-arc bridge-binary-state network. After implementing the BAT, we can have all vectors $X = (x_1, x_2, x_3, x_4, x_5)$, as shown in columns $x_1$, $x_2$, $x_3$, $x_4$, and $x_5$ of Table 2. The details of the other columns and the last row are provided in Section 3.1.

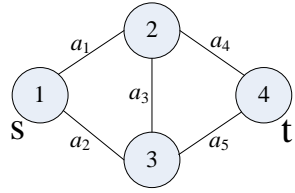

**Figure 2.** The bridge Network

**Table 2.** All vectors obtained from the BAT.

| iteration | $x_1$ | $x_2$ | $x_3$ | $x_4$ | $x_5$ | connect? | $\Pr(x_1)$ | $\Pr(x_2)$ | $\Pr(x_3)$ | $\Pr(x_4)$ | $\Pr(x_5)$ | $\prod_{i=1}^{5}\Pr(x_i)$ |
|---|---|---|---|---|---|---|---|---|---|---|---|---|
| 1 | 0 | 0 | 0 | 0 | 0 | | | | | | | |
| 2 | 1 | 0 | 0 | 0 | 0 | | | | | | | |
| 3 | 0 | 1 | 0 | 0 | 0 | | | | | | | |
| 4 | 1 | 1 | 0 | 0 | 0 | | | | | | | |
| 5 | 0 | 0 | 1 | 0 | 0 | | | | | | | |
| 6 | 1 | 0 | 1 | 0 | 0 | | | | | | | |
| 7 | 0 | 1 | 1 | 0 | 0 | | | | | | | |
| 8 | 1 | 1 | 1 | 0 | 0 | | | | | | | |
| 9 | 0 | 0 | 0 | 1 | 0 | | | | | | | |
| 10 | 1 | 0 | 0 | 1 | 0 | Y | 0.95 | 0.10 | 0.15 | 0.80 | 0.25 | 0.002850 |
| 11 | 0 | 1 | 0 | 1 | 0 | | | | | | | |
| 12 | 1 | 1 | 0 | 1 | 0 | Y | 0.95 | 0.90 | 0.15 | 0.80 | 0.25 | 0.025650 |
| 13 | 0 | 0 | 1 | 1 | 0 | | | | | | | |
| 14 | 1 | 0 | 1 | 1 | 0 | Y | 0.95 | 0.10 | 0.85 | 0.80 | 0.25 | 0.016150 |
| 15 | 0 | 1 | 1 | 1 | 0 | Y | 0.05 | 0.90 | 0.85 | 0.80 | 0.25 | 0.007650 |
| 16 | 1 | 1 | 1 | 1 | 0 | Y | 0.95 | 0.90 | 0.85 | 0.80 | 0.25 | 0.145350 |
| 17 | 0 | 0 | 0 | 0 | 1 | | | | | | | |
| 18 | 1 | 0 | 0 | 0 | 1 | | | | | | | |
| 19 | 0 | 1 | 0 | 0 | 1 | Y | 0.05 | 0.90 | 0.15 | 0.20 | 0.75 | 0.001013 |
| 20 | 1 | 1 | 0 | 0 | 1 | Y | 0.95 | 0.90 | 0.15 | 0.20 | 0.75 | 0.019238 |
| 21 | 0 | 0 | 1 | 0 | 1 | | | | | | | |
| 22 | 1 | 0 | 1 | 0 | 1 | Y | 0.95 | 0.10 | 0.85 | 0.20 | 0.75 | 0.012113 |
| 23 | 0 | 1 | 1 | 0 | 1 | Y | 0.05 | 0.90 | 0.85 | 0.20 | 0.75 | 0.005738 |
| 24 | 1 | 1 | 1 | 0 | 1 | Y | 0.95 | 0.90 | 0.85 | 0.20 | 0.75 | 0.109013 |
| 25 | 0 | 0 | 0 | 1 | 1 | | | | | | | |
| 26 | 1 | 0 | 0 | 1 | 1 | Y | 0.95 | 0.10 | 0.15 | 0.80 | 0.75 | 0.008550 |
| 27 | 0 | 1 | 0 | 1 | 1 | Y | 0.05 | 0.90 | 0.15 | 0.80 | 0.75 | 0.004050 |
| 28 | 1 | 1 | 0 | 1 | 1 | Y | 0.95 | 0.90 | 0.15 | 0.80 | 0.75 | 0.076950 |
| 29 | 0 | 0 | 1 | 1 | 1 | | | | | | | |
| 30 | 1 | 0 | 1 | 1 | 1 | Y | 0.95 | 0.10 | 0.85 | 0.80 | 0.75 | 0.048450 |
| 31 | 0 | 1 | 1 | 1 | 1 | Y | 0.05 | 0.90 | 0.85 | 0.80 | 0.75 | 0.022950 |
| 32 | 1 | 1 | 1 | 1 | 1 | Y | 0.95 | 0.90 | 0.85 | 0.80 | 0.75 | 0.436050 |
| SUM | | | | | | | | | | | | 0.941763 |

To validate the performance of the BAT, it is compared with the breadth-search-first algorithm (BFS) [14], universal generating function method (UGFM) [26], depth-search-first algorithms (DFS) [26], QIE [14], recursive BFS-based SDP (RSDP) [14], and binary-decision diagram (BBD) [27]. Hence, a BAT implemented in the proposed BAT-SSOA3 to calculate the reliability of the general RRAP.

## 2.4 Simplified Swarm Optimization

The first SSO integrated both evolutionary computation and swarm intelligence, which was also proposed by Yeh [20, 23]. Like most AI, all solutions are randomly initialized in the first generation [2, 3, 5, 11, 20, 22, 23]. $x_{i,j}$ is assigned as the $j$th variable of the $i$th solution for $i$ = 1, 2, …, $N_{sol}$ and $j$ = 1, 2, …, $N_{var}$. The basic idea of SSO is to update $x_{i,j}$ for all $i$ and $j$ based on the $i$th variable of the best solution among all the solutions, that is, $G = P_{gBest}$ (called gBest), the best solution in its own evolution history, that is, $P_i$ (called pBest), itself, and a new random feasible value according to the probability $c_g$, $c_p$, $c_w$, and $c_r$, where $c_g + c_p + c_w + c_r = 1$ [20, 23]. The SSO update mechanism is described as follows [20, 23]:

$$x_{i,j} = \begin{cases} g_j & \text{if } \rho_{[0,1]} \in [0, c_g = C_g) \\ p_{i,j} & \text{if } \rho_{[0,1]} \in [c_g, c_g + c_p = C_p) \\ x_{i,j} & \text{if } \rho_{[0,1]} \in [c_g + c_p, c_g + c_p + c_w = C_w) \\ x & \text{if } \rho_{[0,1]} \in [c_g + c_p + c_w, 1] \end{cases} . \tag{7}$$

The SSO pseudocode is presented as follows [20, 23]:

**Algorithm: SSO**

**STEP S0.** Generate $P_i = X_i$ randomly, calculate $F(P_i) = F(X_i)$, let $t$ = 1, and find $gBest$ such that $F(P_i) \le F(P_{gBest})$ for $i$ = 1, 2, …, $N_{sol}$.

**STEP S1.** Let $i$ = 1.

**STEP S2.** Update $X_i$ based on Eq. (7).

**STEP S3.** If $F(P_i) < F(X_i)$, let $P_i = X_i$ and go to STEP S4. Otherwise, go to STEP S5.

**STEP S4.** If $F(P_{gBest}) < F(P_i)$, let $gBest = i$.

**STEP S5.** If $i < N_{sol}$, let $i = i + 1$ and go to STEP S2.

**STEP S6.** If $t < N_{gen}$, let $t = t + 1$ and go to STEP S1. Otherwise, $G = P_{gBest}$ is the final solution.

The flowchart of the SSO is shown in Fig. 3 [20, 23].

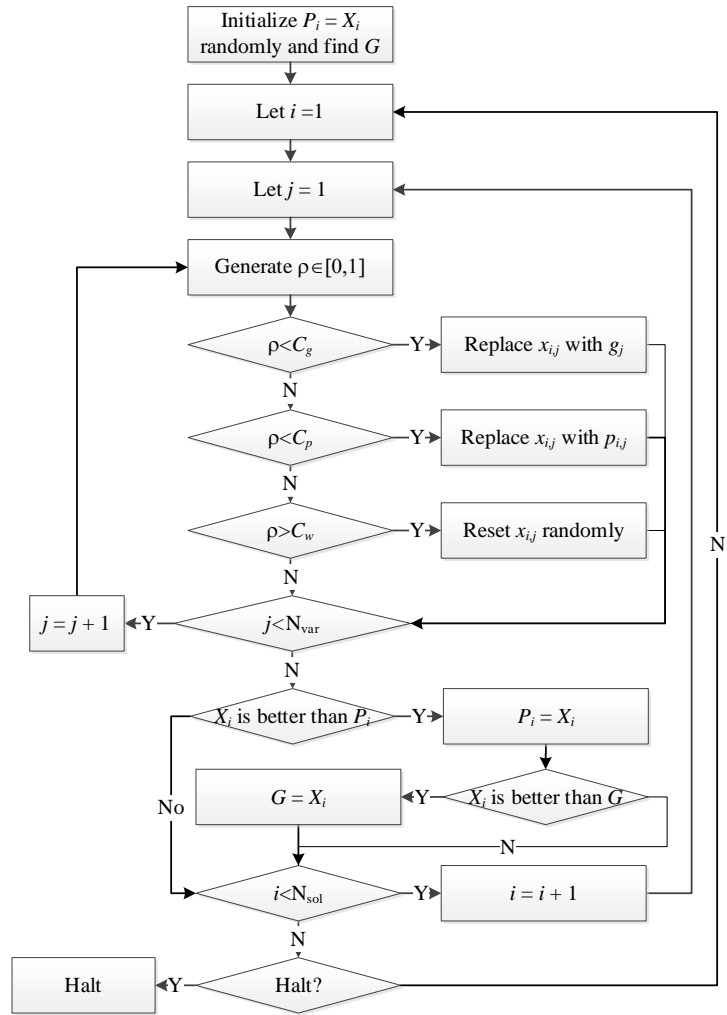

**Figure 3.** SSO flowchart.

SSO has been widely implemented to solve many real-life problems, including disassembly sequencing problems [28], training artificial neural networks (ANNs) [29], high-dimensional multivariable and multimodal numerical continuous benchmark functions [30], parallel-machine scheduling problems [31], dispatch problems [32], RAP/RRAP [1, 2, 3], data mining [33], and allocation and design of RFID networks [34].

Owing to its simplicity and flexibility, SSO has been hybrid with simulated annealing [1], GA [5], PSO [22], ANN [29], bacterial foraging algorithm [35], and artificial bee colony [36]. Furthermore, SSO outperforms PSO, GA, EDA, and ANN in both efficiency and effectiveness (solution quality) based on the computation results [1, 5, 22, 29, 35, 36]. Hence, SSO is adapted in this study to solve the GRRAP.

## 2.5 Penalty Function

There are only three constraints: volume, cost, and weight. If any solution violates any of these constraints is infeasible. To guide solutions to unexplored areas near the boundary of the solution space to obtain an optimal or near-optimal solution, the most common approach is to implement a penalty function for these problems with few constraints [2, 3, 5, 6, 21, 22]. The most common penalty function is as follows:

$$R_p(\mathbf{n},\ \mathbf{r}) = R_s(\mathbf{n},\ \mathbf{r}) \cdot \left( Min\left\{ \left[ \frac{V_{ub}}{g_v(\mathbf{n},\ \mathbf{r})} \right], \left[ \frac{C_{ub}}{g_c(\mathbf{n},\ \mathbf{r})} \right], \left[ \frac{W_{ub}}{g_w(\mathbf{n},\ \mathbf{r})} \right], 1 \right\} \right)^3. \tag{8}$$

The system/network reliability $R_s(\mathbf{n},\ \mathbf{r})$ of any infeasible solution is replaced by the penalty function $R_p(\mathbf{n},\ \mathbf{r})$ in the proposed BAT-SSOA3.

**Example 2.** Assume that $R_s(X) = .9316$, $g_v(X) = 199$, $g_c(X) = 174.9$, and $g_w(X) = 30.9$ and let $V_{UB} = 200$, $C_{UB} = 175$, and $W_{UB} = 30$. Thus,

$$R_p(X) = R_s(X) \cdot \left( Min\left\{ \left[ \frac{V_{ub}}{g_v(X)} \right], \left[ \frac{C_{ub}}{g_c(X)} \right], \left[ \frac{W_{ub}}{g_w(X)} \right], 1 \right\} \right)^3$$

$$= 0.8525. \tag{9}$$

## 3. THE PROPOSED BAT-SSOA3

There are four major innovations in the proposed **BAT-SSOA3**.

1. Select parameters and terms in Eq. (7) based on the training using small samples for all test problems in the proposed SS3OA;

2. Calculate the reliability using the connected vectors obtained from the BAT;

3. Update $G$ using the boundary update.

4. Add the redundancy variable $n_i$ and reliability variable $r_i$ as variable $x_i = (n_i + r_i)$ to reduce the update time.

These four parts are discussed and proposed in this section to improve the SSO for mixed-integer nonlinear programming problems.

### 3.1 BAT-Connected Vectors in Calculating Network Reliability

All vectors in the solution space can be obtained after implementing the BAT [14, 25, 26, 27], as discussed in Section 2.3. To ensure optimal use of computer memory, BAT uses only one vector from the beginning to the end of the entire process. However, this advantage is modified in which all connected vectors are saved to calculate the reliability of the GRRAP without finding all connected vectors every time the objective function is calculated.

A vector $X$ is considered connected if its corresponding $G(X)$ is connected, that is, there is a path from the input node to the output node in $G(X)$. An example is assigning nodes 1 and 4 as the input and output nodes, respectively. Fig. 4(1) and 4(2) show $G(X)$ and $G(Y)$, where $X = (0, 0, 1, 0, 1)$ and $Y = (0, 0, 1, 0, 1)$. In addition, $X$ is disconnected, and $Y$ is connected.

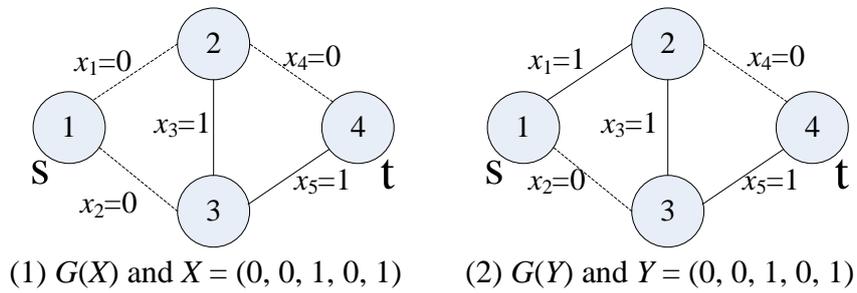

(1) $G(X)$ and $X = (0, 0, 1, 0, 1)$     (2) $G(Y)$ and $Y = (0, 0, 1, 0, 1)$
**Figure 4.** Examples of connected vectors.

In the proposed BAT, all connected vectors, for example, vectors 10, 12, 14 – 16, 19, 20, 22 – 24, 26 – 28, 30 – 32, which are marked "Y" in the column titled "connect?" in Table 2, are stored to calculate the system/network reliability for each solution at each generation. $\Pr(x_i)$ denotes the probability when the state of $a_i = x_i$, and $\prod_{i=1}^{5} \Pr(x_i)$ is the probability of the connected vector $X = (x_1, x_2, x_3, x_4, x_5)$.

For example, using the data listed in the second last row of Table 2, $\Pr(X) = 0.95 \times 0.90 \times 0.85 \times 0.80 \times 0.75 = 0.436050$, where $\Pr(x_1) = 0.95$, $\Pr(x_2) = 0.90$, $\Pr(x_3) = 0.85$, $\Pr(x_4) = 0.80$, and $\Pr(x_5) = 0.75$. Finally, the sum of the probabilities of all complete vectors, that is, $\sum_{X} \prod_{i=1}^{5} \Pr(x_i) = 0.941763$ for all connected vectors $X$, is the reliability.

In addition, prior to calculating the reliability, each test problem needs to run the BAT once only to have all the connected vectors.

## 3.2 Small-Sampling Tri-Objective Orthogonal Array

Setting and tuning parameters systematically and efficiently without exploring all possible combinations are important for most AI. The proposed small-sampling tri-objective orthogonal array (SS3OA) is based on an orthogonal array (OA) and small sampling and has three goals:

1. Tune parameters (i.e., $C_g$, $C_p$, and $C_w$) in SSO systematically.

2. Decide which items (i.e., items 1–4 in Eq. (7)) can be removed to efficiently obtain a better solution quality.

3. Avoid misleading and overfitting in using the OA.

An OA is an array in the design of experiments [1, 37], and it is used to study several factors simultaneously to determine the best combination of factor levels efficiently without trying all possible combinations. It arranges factors, such as $c_g$, $c_p$, and $c_w$, in the proposed SS3OA in columns and trys in rows, such that each column is orthogonal to the other columns. $L_a(b^c)$ is a standard OA notation with $a$ trys, $b$ levels, and $c$ factors. To test whether it is better to have higher, medium, and low levels of the corresponding factors, each parameter has three levels: 1, 2, and 3, which correspond to large, medium, and low, respectively.

$L_9(3^3)$ shown in Table 3 is used in the proposed SS3OA, and the related level values are shown in Table 4.

**Table 3.** Example $L_9(3^3)$ and the levels of each factor.

| Try | Factor 1 | Factor 2 | Factor 3 |
|-----|----------|----------|----------|
| 1 | 1 | 1 | 1 |
| 2 | 1 | 2 | 2 |
| 3 | 1 | 3 | 3 |
| 4 | 2 | 1 | 2 |
| 5 | 2 | 2 | 3 |
| 6 | 2 | 3 | 1 |
| 7 | 3 | 1 | 3 |
| 8 | 3 | 2 | 1 |
| 9 | 3 | 3 | 2 |

**Table 4.** Three levels and the related probabilities.

| | $c_g$ | $c_p$ | $c_w$ |
|---|-------|-------|-------|
| 1 | 0.6 | 0.30 | 0.30 |
| 2 | 0.4 | 0.25 | 0.20 |
| 3 | 0.2 | 0 | 0 |

In addition, to achieve our second goal, we tested removing the second, third, and/or fourth items in Eq. (7) and set some levels to zero. For example, $c_p$ and $c_w$ at level 3 denote that the new updated variable is either from $G$ (the first item in Eq. (7)) or from a newly generated feasible random number (the fourth item in Eq. (7)).

For example, after listing $C_g = c_g$, $C_p = c_g + c_p$, and $C_w = C_p + c_w$ in Tables 3 and 4, Table 5 was obtained. In Table 5, Trys 0 and 1 had $C_w > 1.0$, which indicates that $c_r = 0$, that is, no random feasible value is generated to replace the related variable. Trys 2, 5, and 8 had $C_g = C_p$, which indicates that the role of *pBest* is completely removed, that is, Item 2 is discarded in Eq. (1). Runs 2, 4, and 6 had $C_p = C_w$, which indicates that the third item in Eq. (1) is removed.

**Table 5.** Three levels and the related probabilities.

| Try | Factor 1 | Factor 2 | Factor 3 | $c_g$ | $c_p$ | $c_w$ | $C_g$ | $C_p$ | $C_w$ |
|-----|----------|----------|----------|-------|-------|-------|-------|-------|-------|
| 0 | 1 | 1 | 1 | 0.6 | 0.30 | 0.30 | 0.60 | 0.90 | 1.20 |
| 1 | 1 | 2 | 2 | 0.6 | 0.25 | 0.20 | 0.60 | 0.85 | 1.05 |
| 2 | 1 | 3 | 3 | 0.6 | 0 | 0 | 0.60 | 0.60 | 0.60 |
| 3 | 2 | 1 | 2 | 0.4 | 0.30 | 0.20 | 0.40 | 0.70 | 0.90 |
| 4 | 2 | 2 | 3 | 0.4 | 0.25 | 0 | 0.40 | 0.65 | 0.65 |
| 5 | 2 | 3 | 1 | 0.4 | 0 | 0.30 | 0.40 | 0.40 | 0.70 |
| 6 | 3 | 1 | 3 | 0.2 | 0.30 | 0 | 0.20 | 0.50 | 0.50 |
| 7 | 3 | 2 | 1 | 0.2 | 0.25 | 0.30 | 0.20 | 0.45 | 0.75 |
| 8 | 3 | 3 | 2 | 0.2 | 0 | 0.20 | 0.20 | 0.20 | 0.40 |

After performing all nine Trys, the next step for the SS3OA is to average these fitness values based on the levels of each factor. For example, in Table 6, the first three values below Try 8, that is, 0.9997044, 0.9999069, and 0.9957876, are the averages of the fitness values based on Levels 1, 2, and 3. That is, 0.9997044 = (0.9997173 + 0.9995158 + 0.9998801)/3 as 0.9997173, 0.9995158, and 0.9998801 are obtained based on Level 1 of Factor 1 (i.e., $c_g$); 0.9999069 = (0.9999575 + 0.9998977 + 0.9998655)/3 as 0.9999575, 0.9998977, and 0.9998655 are obtained based on Level 2 of Factor 1; 0.9957876 = (0.9991318 + 0.9998328 + 0.9883984)/3 as 0.9991318, 0.9998328, and 0.9883984 are obtained based on Level 3 of Factor 1. In addition, 0.9999069 is the best among the three averages. Hence, there is a new Try, that is, Try 9, whose level for Factor 1 is Level 2, that is, $c_g = 0.4$. Similarly, the levels for Factors 2 and 3 are levels 2 and 1 in Try 9, respectively.

**Table 6.** Average fitness values based on levels.

| | Try | Fitness |
|---|---|---|
| | 0 | 0.9997173 |
| | 1 | 0.9995158 |
| | 2 | 0.9998801 |
| | 3 | 0.9999575 |
| | 4 | 0.9998977 |
| | 5 | 0.9998655 |
| | 6 | 0.9991318 |
| | 7 | 0.9998328 |
| | 8 | 0.9883984 |
| $c_g$ | 1 | 0.9997044 |
| | 2 | **0.9999069** |
| | 3 | 0.9957876 |
| $c_p$ | 1 | 0.9996022 |
| | 2 | **0.9997488** |
| | 3 | 0.9960480 |
| $c_w$ | 1 | **0.9998052** |
| | 2 | 0.9959572 |
| | 3 | 0.9996365 |

To achieve our third goal for the proposed SS3OA, traditional parameters could either be trained from a problem and implemented to test other problems, or trained from all test problems and tested on all problems again. The former could result in misleading information, whereas the latter could result in overfitting. $N_{run}$ is assigned as the number of runs for each test problem. To overcome both issues in the proposed SS3OA, we only train five independent runs for all problems.

**3.3 Boundary Update**

Based on experience and from literature, SSO is powerful in global searches, but may be weak in local search. The boundary update is derived by differentiating the cost function, as shown in Eq. (5), and was first proposed in [2, 3]. Owing to their effectiveness, some SSO algorithms have already implemented boundary updates to improve the local search of SSO [2, 3]. Hence, the boundary update is also adapted in the proposed SS3OA; however, it limited the update of one variable in $G$ at each generation as it required more time to implement Eq. (5).

The pseudocode for the proposed boundary update to update the current $G$ is as follows:

**Algorithm: Boundary_Update**

**STEP U0.** Let the integer remainder of gen/Nvar be $i$ and $G = (g_1, g_2, \ldots, g_{Nvar})$.

**STEP U1.** Update $g_i$ to $g^*$ based on the boundary condition listed in Eq. (1).

**STEP U2.** If new $G$, say $G^*$, is improved, let $G = G^*$; otherwise, keep $G$ and discard the update.

Note that boundary updates have been used in other (continuous) SSO variants, as in [2, 3]. However, they are first used in the traditional (discrete) SSO, which update solutions simply based on Eq. (7) in the proposed BAT-SSOA3.

## 3.4 Affixation Solution Structure

Both the redundancy variables (integer variables) and reliability variables (floating-point variables) are included in traditional RRAP and the proposed GRRAP. Thus, we need two different update mechanisms, that is, one for integer variables and one for floating-point variables, to update different types of variables. To simplify the update mechanism to reduce runtime, the proposed BAT-SSOA3 combines two different types of variables. For example, if $n_i$ and $r_i$ are from the same subsystem, they could be added, that is, $(n_i + r_i)$, to create one floating-point variable. For instance, in Example 1 where $\mathbf{n} = (3, 2, 2, 3, 3)$ and $\mathbf{r} = (0.77946645, 0.87173278, 0.90284951, 0.71148780, 0.78781644)$, they can be combined to form $X = \mathbf{n} + \mathbf{r} = (3.77946645, 2.87173278, 2.90284951, 3.71148780, 3.78781644)$.

## 3.5 Pseudocode of the Proposed BAT-SSOA3

The detailed procedure for the proposed BAT-SSOA3 is as follows:

**STEP 0.**   Based on SS3OA, run the proposed BAT-SSOA3 five times for each test problem, and each run has nine different parameter settings, as shown in Table 5. Find the combination that has the highest average for the last $G$ in all test problems and in all five runs with nine settings.

**STEP 1.**   Run BAT and store all connected vectors.

**STEP 2.**   Generate $P_i = X_i = (\mathbf{r}_i + \mathbf{n}_i)$ randomly, calculate $R_p(P_i) = R_p(X_i)$ based on the penalty function and connected vectors obtained from BAT, let $t = 1$, and find $G$ such that $R_p(P_i) \leq R_p(G)$ for $i = 1, 2, \ldots, N_{sol}$.

**STEP 3.**   Let $i = 1$.

**STEP 4.**   Update $X_i$ based on Eq. (7) if $X_i \neq G$. Otherwise, update $X_i$ based on the proposed boundary update.

**STEP 5.**    If $R_p(P_i) < R_p(X_i)$, let $P_i = X_i$. Otherwise, go to STEP 7.

**STEP 6.**    If $R_p(G) < R_p(P_i)$, let $G = P_i$.

**STEP 7.**    If $i < N_{sol}$, let $i = i + 1$ and go to STEP S2.

**STEP 8.**    If $t < N_{gen}$, let $t = t + 1$ and go to STEP S1. Otherwise, $P_{gBest}$ is the final solution.

It is trivial that the SSO in the proposed BAT-SSOA3 is similar to the traditional SSO, as in STEPs 2–8, except that

1.  SS3OA is implemented to select the values of $C_g$, $C_p$, and $C_w$, and determine which terms in Eq. (7) must be discarded or kept as in STEP 0.

2.  The connected vectors determined and stored in the BAT are used to calculate $R_p(X_i)$ in STEP 1.

3.  The penalty function is adopted to calculate $R_p(X_i)$ in STEPs 2 and 4.

4.  $G$ is updated using the proposed boundary update.

5.  All solutions, including $X_i$, $P_i$, and $G$, are constructed using the affixation solution structure, as discussed in Section 3.4.

## 4. Numerical examples

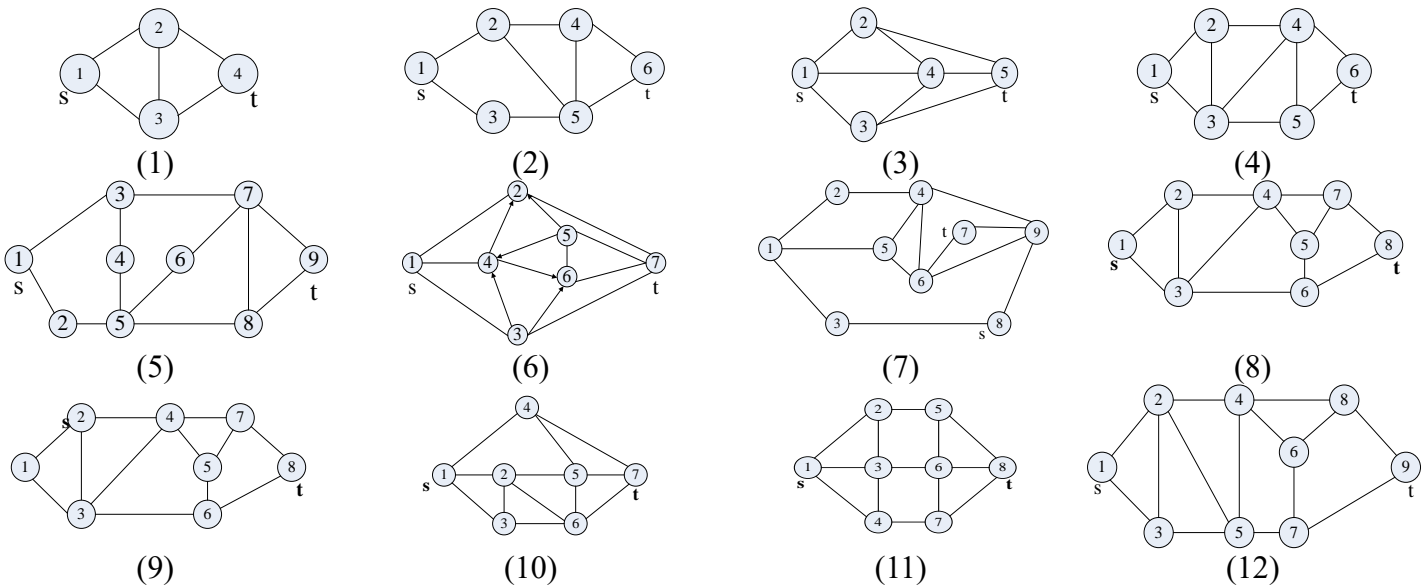

**Figure 5.** Twelve test problems.

To evaluate the algorithm performance and solution quality of the proposed BAT-SSOA3 for the proposed general RRAPs, two experiments, Ex1 and Ex2, are included in this study. In addition, 12 test problems, shown in Fig. 5, are used to validate the performance of the BAT-SSOA3 by comparing it with other state-of-the-art algorithms. Note that the 12 problems are well-known binary-state networks in the network research area [14, 25, 26, 27].

## 4.1 Test Environment

Each subsystem in each test problem of the two experiments can have a maximum of 10 identical components in parallel. The corresponding constants of the $i$th arc in each network, which are $\alpha_i \cdot 10^5$, $\beta_i$, $w_i v_i^2$, and $w_i$ in Eq. (2)–(4), are similar to those of the $i$th subsystem in Table 1, where

$$\iota = \begin{cases} 5 & \text{if the remainder of } i/5 \text{ is zero} \\ \text{the reminder of } i/5 & \text{otherwise} \end{cases}. \tag{10}$$

For example, in the 6$^{\text{th}}$ arc of each network for all test problems, $\alpha_i \cdot 10^5 = 2.330$, $\beta_i = 1.5$, $w_i v_i^2 = 1$, and $w_i = 7$, which are all copied from the 1$^{\text{st}}$ subsystem in Table 1, as the remainder of 6/5 is 1.

The three upper bounds of constraints $V_{\text{ub}}$, $C_{\text{ub}}$, and $W_{\text{ub}}$ are also based on Table 1, such that

$$\bullet_{\text{ub,new}} = N_{\text{var}} \times \bullet_{\text{ub,old}}/5. \tag{11}$$

For example, in a network with nine arcs (Fig. 5(4)), $V_{\text{ub}} = 9 \times 110/5 = 198$, $C_{\text{ub}} = 9 \times 175/5 = 315$, and $W_{\text{ub}} = 9 \times 200/5 = 360$ as $V_{\text{ub}} = 110$, $C_{\text{ub}} = 175$, and $W_{\text{ub}} = 200$ in Table 1.

Ex 1 had ten settings for the proposed BAT-SSOA3, that is, there are ten different BAT-SSOA3, whereas Ex 2 had six algorithms. Hence, there were 10 + 6 = 15 algorithms, including the proposed BAT-SSOA3 in Ex 1 and Ex 2.

All the algorithms in both Ex 1 and Ex 2 are programmed in C$^{++}$ using Dev C++5.11, carried out on an Intel Core i9-9900K @3.60 GHz PC with 48 GB memory Windows 10 (runtime unit is in CPU seconds). In addition, for a fair comparison, in each test, $N_{\text{run}} = 100$ and $N_{\text{gen}} = 1000$, which is also the stopping criterion, $N_{\text{sol}} = 100$, and all algorithms had the same calculation number of fitness (reliability).

Moreover, to easily observe the convergent trend, we separated $N_{\text{gen}} = 1000$ into four equal-size generations, that is, 1–250, 251–500, 501–750, and 751–1000, and recorded their corresponding results.

Moreover, the 250[th], 500[th], 750[th], and 1000[th] generations are called stages 0, 1, 2, and 3, respectively, in the rest of the tables.

## 4.2 Ex 1

Ex 1 mainly tests and verifies the benefit of the proposed SS3OA, which not only systematically selects the values of $C_g$, $C_p$, and $C_w$ but also efficiently decides which terms in Eq. (7) should be retained or discarded. The BAT-SSOA3 with the best parameter settings and terms in Eq. (7) from SS3OA, of Ex 1 for all test problems are adopted in Ex 2 for comparison with other algorithms.

As discussed in Section 3.2, there are only five runs for each test problem in the proposed SS3OA, where each problem had (9+1) different parameter settings, as shown in Tables 7–9.

Note that in these tables, stages 0, 1, 2, and 3 represent the 250[th], 500[th], 750[th], and 1000[th] generations, respectively. The bold values are the best among the nine trials. The number in bold indicates the best result among the nine tests of the same stage for the same test problem. For example, based on Table 7, 0.999781 is the best among the nine tests in test problem 1 at stage 0.

Based on the results listed in Tables 7 and 8, the average of the best final $G$ or the mean averages obtained from BAT-SSOA3 are from Try 3, that is, $C_g = 0.4$, $C_p = 0.7$, and $C_w = 0.9$. In addition, Try 3 had a greater number of averages for the best final $G$ and had the best mean averages listed in Tables 7–9. In contrast, these tries with $c_r = 0$, that is, Tries 0 and 1, are already in local optima and convergence in the first stage. The largest $c_r$ indicates the worst solution quality for all aspects; for example, $c_r = 0.6$ for Try 8 (the worst) and $c_r = 0.5$ for Try 6 (the second worst). Hence, it is not appropriate to discard the fourth item in Eq. (7), and a larger value for $c_r$. The reason is that a smaller/zero $c_r$ indicates a smaller/zero chance to escape the local trap after the solution is trapped; larger $c_r$ indicates a larger divergence and leaves the optimum far away.

Based on Table 9, $C_g = 0.40$ was preferable, for all stages, as an average $c_g$ indicates high reliability. Hence, a middle value of $C_g$ indicates an improved probability of obtaining a better solution. However, for both $c_p$ and $c_w$, the best levels were not fixed values. For example, $c_p$ was better when levels were 2, 2, 1, and 1 for stages 0–3, respectively, and $c_w$ was better when levels were 1, 3, 3, and 3 for stages 0–3,

respectively

As discussed in Section 3.2, to decide whether to adapt and implement $C_g = 0.4$, $C_p = 0.7$, and $C_w = 0.9$ in Ex 2 for the proposed BAT-SSOA3, the average based on the levels of each factor, which are summarized in Table 9, must be determined. Thereafter, all settings are completed to determine the best levels for $c_g$, $c_p$, and $c_w$ for obtaining the best solutions listed in Table 10.

**Table 7.** Average of the best $G$ for each test problem in each stage for five runs.

| stage | Try | 1 | 2 | 3 | 4 | 5 | 6 | 8 | 9 | 10 | 11 | 12 | 15 | average |
|---|---|---|---|---|---|---|---|---|---|---|---|---|---|---|
| 0 | 0 | 0.99911 | 0.999527 | 1 | 0.999629 | 0.999414 | 0.999983 | 0.999936 | 0.999731 | 0.999539 | 0.999997 | 0.999999 | 0.999745 | 0.999717 |
| | 1 | 0.999599 | 0.999811 | 0.999999 | 0.999494 | 0.999553 | 0.999995 | 0.999913 | 0.999753 | 0.996897 | 0.999987 | 0.999959 | 0.999228 | 0.999516 |
| | 2 | **0.999781** | 0.99991 | **1** | 0.999939 | 0.99986 | 0.999999 | 1 | 0.999923 | 0.999302 | 0.999998 | 0.999999 | 0.999852 | 0.99988 |
| | 3 | 0.999701 | **0.999957** | **1** | **0.999957** | **0.999982** | 1 | 1 | **0.999977** | **0.999969** | 1 | 1 | **0.999946** | **0.999958** |
| | 4 | 0.99978 | 0.999948 | 1 | 0.999933 | 0.999911 | 0.999998 | 0.999997 | 0.999889 | 0.999582 | 0.999999 | 0.999996 | 0.999739 | 0.999898 |
| | 5 | 0.999773 | 0.999918 | 1 | 0.999843 | 0.999945 | 0.999997 | 0.999996 | 0.999862 | 0.99946 | 0.999998 | 0.999995 | 0.999599 | 0.999865 |
| | 6 | 0.999735 | 0.999798 | 1 | 0.999523 | 0.998638 | 0.999976 | 0.99981 | 0.999441 | 0.996502 | 0.999967 | 0.99983 | 0.9621 | 0.999132 |
| | 7 | 0.99978 | 0.999915 | **1** | 0.999942 | 0.999753 | 0.999996 | 0.999993 | 0.99978 | 0.999513 | 0.999999 | 0.999996 | 0.999325 | 0.999833 |
| | 8 | 0.999588 | 0.999541 | 0.99999 | 0.998404 | 0.986951 | 0.999304 | 0.992515 | 0.999157 | 0.913531 | 0.999967 | 0.999156 | 0.972676 | 0.988398 |
| | 9 | 0.999767 | **0.999963** | **1** | 0.999984 | | 1 | | 0.999979 | 0.999962 | 1 | 1 | **0.999959** | **0.999965** |
| 1 | 0 | 0.99911 | 0.999527 | 1 | 0.999629 | 0.999414 | 0.999983 | 0.999936 | 0.999731 | 0.999539 | 0.999997 | 0.999999 | 0.999745 | 0.999717 |
| | 1 | 0.999599 | 0.999811 | 0.999999 | 0.999494 | 0.999553 | 0.999995 | 0.999913 | 0.999753 | 0.996897 | 0.999987 | 0.999959 | 0.999228 | 0.999516 |
| | 2 | **0.999781** | 0.999936 | 1 | 0.999952 | 0.999967 | 0.999999 | 1 | 0.999956 | 0.999675 | 0.999999 | 0.999999 | 0.999884 | 0.999939 |
| | 3 | 0.999715 | **0.999967** | **1** | 0.999959 | **0.999983** | 1 | 1 | **0.999982** | **0.999974** | 1 | 1 | **0.999951** | **0.999961** |
| | 4 | 0.99978 | 0.999962 | 1 | **0.999965** | 0.999971 | 0.999999 | 0.999998 | 0.99995 | 0.999826 | 0.999999 | 0.999997 | 0.999879 | 0.999944 |
| | 5 | 0.999774 | 0.999931 | 1 | 0.99991 | 0.999947 | 0.999997 | 0.999997 | 0.999934 | 0.999521 | 0.999998 | 0.999998 | 0.999775 | 0.999901 |
| | 6 | 0.999748 | 0.999922 | 1 | 0.999857 | 0.99923 | 0.999985 | 0.999836 | 0.999441 | 0.998901 | 0.999996 | 0.999983 | 0.999333 | 0.999686 |
| | 7 | 0.99978 | 0.999963 | 1 | 0.999957 | 0.999864 | 0.999999 | 0.999999 | 0.999915 | 0.999657 | 0.999999 | 0.999997 | 0.999813 | 0.999912 |
| | 8 | 0.999708 | 0.999541 | 0.999996 | 0.998692 | 0.99987 | 0.999879 | 0.998737 | 0.999917 | 0.973943 | 0.999972 | 0.999156 | | 0.995232 |
| | 9 | 0.999767 | **0.999967** | **1** | 0.999958 | **0.999985** | 1 | 1 | 0.99984 | 0.999973 | 1 | 1 | 0.99962 | **0.999966** |
| 2 | 0 | 0.99911 | 0.999527 | 1 | 0.999629 | 0.999414 | 0.999983 | 0.999936 | 0.999731 | 0.999539 | 0.999997 | 0.999999 | 0.999745 | 0.999717 |
| | 1 | 0.999599 | 0.999811 | 0.999999 | 0.999494 | 0.999553 | 0.999995 | 0.999913 | 0.999753 | 0.996897 | 0.999987 | 0.999959 | 0.999228 | 0.999516 |
| | 2 | **0.999781** | 0.999957 | 1 | 0.999952 | 0.999967 | 0.999999 | 1 | 0.999968 | 0.999766 | 0.999999 | 0.999999 | 0.999884 | 0.999939 |
| | 3 | 0.999757 | **0.99997** | **1** | 0.999962 | **0.999985** | 1 | 1 | **0.999985** | **0.999976** | 1 | 1 | **0.999959** | **0.999966** |
| | 4 | 0.99978 | 0.999967 | 1 | **0.999965** | 0.99997 | 1 | 0.999999 | 0.999961 | 0.999873 | 1 | 0.999999 | 0.999885 | 0.999951 |
| | 5 | 0.99978 | 0.99994 | 1 | 0.999922 | 0.99997 | 0.999999 | 0.999999 | 0.999934 | 0.999729 | 0.999999 | 0.999998 | 0.999775 | 0.999921 |
| | 6 | 0.99976 | 0.999922 | 1 | 0.999897 | 0.999439 | 0.999993 | 0.999836 | 0.999583 | 0.998901 | 0.999996 | 0.999983 | 0.999333 | 0.99972 |
| | 7 | 0.99978 | 0.999963 | 1 | 0.999961 | 1 | 0.999999 | 0.999999 | 0.999915 | 0.999789 | 0.999999 | 0.999998 | 0.999813 | 0.999936 |
| | 8 | 0.999753 | 0.999541 | 0.999998 | 0.999045 | 0.986951 | 0.999879 | 0.998737 | 0.999375 | 0.990369 | 0.999972 | 0.999876 | 0.990021 | 0.9696 |
| | **9** | **0.99978** | **0.99997** | **1** | **0.999966** | **0.999985** | 1 | 1 | **0.999985** | 0.999975 | 1 | 1 | **0.999968** | **0.999969** |
| 3 | 0 | 0.99911 | 0.999527 | 1 | 0.999629 | 0.999414 | 0.999983 | 0.999936 | 0.999731 | 0.999539 | 0.999997 | 0.999999 | 0.999745 | 0.999717 |
| | 1 | 0.999599 | 0.999811 | 0.999999 | 0.999494 | 0.999553 | 0.999995 | 0.999913 | 0.999753 | 0.996897 | 0.999987 | 0.999959 | 0.999228 | 0.999516 |
| | 2 | 0.999781 | 0.999963 | 1 | 0.99996 | 0.999968 | 0.999999 | 1 | 0.999968 | 0.999893 | 1 | 1 | 0.999884 | 0.999951 |
| | 3 | 0.999758 | 0.99997 | 1 | **0.999966** | **0.999985** | 1 | 1 | **0.999986** | **0.999976** | 1 | 1 | 0.99996 | **0.999967** |
| | 4 | **0.999782** | **0.99997** | **1** | 0.999965 | 0.99979 | 1 | 0.999999 | 0.999978 | 0.999906 | 1 | 0.999999 | 0.999895 | 0.999956 |
| | 5 | 0.99978 | 0.99994 | 1 | 0.999924 | 0.99997 | 0.999999 | 1 | 0.999934 | 0.999729 | 0.999999 | 0.999998 | 0.999857 | 0.999928 |
| | 6 | 0.999767 | 0.999924 | 1 | 0.999907 | 0.999996 | 0.999996 | 0.999989 | 0.999624 | 0.99714 | 0.999996 | 0.999984 | 0.999333 | 0.999835 |
| | 7 | 0.99978 | 0.999963 | 1 | **0.999967** | 0.999912 | 0.999999 | 0.999999 | 0.99996 | 0.999806 | 1 | 0.999999 | 0.999883 | 0.999939 |
| | 8 | 0.999776 | 0.999541 | 0.999998 | 0.999045 | 0.989012 | 0.999879 | 0.998737 | 0.999375 | 0.990369 | 0.999972 | 0.999876 | 0.997185 | 0.99773 |
| | 9 | 0.999781 | **0.99997** | **1** | 0.999966 | **0.999985** | 1 | 1 | 0.999985 | 0.999975 | **1** | **1** | **0.999968** | **0.999969** |

**Table 8.** Average of the average $G$ for each test problem in each stage for five runs.

| stage | Try | 1 | 2 | 3 | 4 | 5 | 6 | 8 | 9 | 10 | 11 | 12 | 15 | average |
|---|---|---|---|---|---|---|---|---|---|---|---|---|---|---|
| 0 | 0 | 0.997963 | 0.999078 | 0.982713 | 0.999367 | 0.998115 | 0.999967 | 0.990083 | 0.994797 | 0.987132 | 0.996593 | 0.998744 | 0.987148 | 0.994308 |
| | 1 | 0.998412 | 0.998405 | 0.999967 | 0.994694 | 0.989406 | 0.999908 | 0.99938 | 0.999783 | 0.992692 | 0.996846 | 0.99965 | 0.996275 | 0.996955 |
| | 2 | 0.999742 | 0.999827 | 1 | 0.999842 | 0.99962 | 0.999997 | 0.999939 | 0.999757 | 0.99869 | 0.999996 | 0.99999 | 0.999066 | 0.999705 |
| | 3 | 0.999631 | **0.999899** | **1** | **0.999925** | **0.999951** | **1** | **0.999998** | **0.999966** | **0.999876** | **1** | **1** | **0.999931** | **0.999931** |
| | 4 | 0.999732 | 0.999874 | 1 | 0.999831 | 0.999782 | 0.999997 | 0.999986 | 0.999801 | 0.999301 | 0.999997 | 0.999991 | 0.999246 | 0.999795 |
| | 5 | **0.99975** | 0.999809 | 0.999999 | 0.999742 | 0.99955 | 0.999986 | 0.999952 | 0.999794 | 0.998519 | 0.999995 | 0.999965 | 0.999271 | 0.999694 |
| | 6 | 0.999649 | 0.999673 | 0.999997 | 0.999327 | 0.996441 | 0.999891 | 0.997705 | 0.998369 | 0.991468 | 0.99994 | 0.999833 | 0.992929 | 0.997935 |
| | 7 | 0.999682 | 0.999803 | 1 | 0.99988 | 0.999381 | 0.999988 | 0.999977 | 0.999391 | 0.998756 | 0.999994 | 0.999986 | 0.998755 | 0.999633 |
| | 8 | 0.999322 | 0.9977 | 0.999911 | 0.996349 | 0.904061 | 0.830335 | 0.881246 | 0.94601 | 0.679505 | 0.974904 | 0.96735 | 0.767328 | 0.913176 |
| | 9 | 0.999449 | **0.999911** | **1** | 0.999837 | **0.999963** | **1** | **0.999998** | **0.999969** | **0.999936** | **1** | **1** | **0.999945** | 0.999917 |
| 1 | 0 | 0.997963 | 0.999078 | 0.982713 | 0.999367 | 0.998115 | 0.999967 | 0.990083 | 0.994797 | 0.987132 | 0.996593 | 0.998744 | 0.987148 | 0.994308 |
| | 1 | 0.998412 | 0.998405 | 0.999967 | 0.994694 | 0.989406 | 0.999908 | 0.99938 | 0.99783 | 0.992692 | 0.996846 | 0.99965 | 0.996275 | 0.996955 |
| | 2 | **0.999756** | 0.999859 | 1 | 0.999918 | 0.999763 | 0.999999 | 0.999994 | 0.999859 | 0.999465 | 0.999998 | 0.999997 | 0.999483 | 0.999841 |
| | 3 | 0.999645 | 0.999906 | 1 | **0.99995** | **0.999962** | **1** | **0.999998** | **0.999974** | **0.999906** | **1** | **1** | **0.999942** | **0.99994** |
| | 4 | 0.999735 | 0.999906 | 1 | 0.999902 | 0.999913 | 0.999999 | 0.999999 | 0.999901 | 0.999706 | 0.999999 | 0.999994 | 0.999652 | 0.999892 |
| | 5 | 0.999755 | 0.999845 | 1 | 0.999877 | 0.999806 | 0.999993 | 0.999986 | 0.999875 | 0.999053 | 0.999996 | 0.999984 | 0.999615 | 0.999815 |
| | 6 | 0.999694 | 0.999789 | 0.999999 | 0.999701 | 0.998647 | 0.99971 | 0.99878 | 0.998914 | 0.996751 | 0.999975 | 0.999953 | 0.995509 | 0.998974 |
| | 7 | 0.999745 | **0.999931** | 1 | 0.999927 | 0.999799 | 0.999999 | 0.999997 | 0.999799 | 0.999511 | 0.999998 | 0.999993 | 0.999572 | 0.999861 |
| | 8 | 0.999566 | 0.99836 | 0.999988 | 0.998144 | 0.972585 | 0.999152 | 0.97471 | 0.992516 | 0.885962 | 0.99963 | 0.997258 | 0.930947 | 0.979068 |
| | 9 | 0.999536 | **0.999942** | **1** | 0.999933 | **0.999969** | **1** | **0.999999** | **0.999978** | 0.999952 | **1** | **1** | **0.999953** | 0.999938 |
| 2 | 0 | 0.997963 | 0.999078 | 0.982713 | 0.999367 | 0.998115 | 0.999967 | 0.990083 | 0.994797 | 0.987132 | 0.996593 | 0.998744 | 0.987148 | 0.994308 |
| | 1 | 0.998412 | 0.998405 | 0.999967 | 0.994694 | 0.989406 | 0.999908 | 0.99938 | 0.99783 | 0.992692 | 0.996846 | 0.99965 | 0.996275 | 0.996955 |
| | 2 | **0.999764** | 0.99987 | 1 | 0.999931 | 0.999777 | 0.999999 | 0.999995 | 0.999912 | 0.999665 | 0.999999 | 0.999998 | 0.999657 | 0.999881 |
| | 3 | 0.999673 | **0.999943** | **1** | **0.999957** | **0.999964** | **1** | **0.999998** | **0.999979** | **0.999923** | **1** | **1** | **0.999945** | **0.999948** |
| | 4 | 0.999735 | 0.999939 | 1 | 0.999938 | 0.999938 | 0.999999 | 0.999995 | 0.999929 | 0.999767 | 0.999999 | 0.999995 | 0.999682 | 0.99991 |
| | 5 | 0.999757 | 0.999882 | 1 | 0.999896 | 0.999834 | 0.999997 | 0.999995 | 0.999881 | 0.999289 | 0.999998 | 0.999989 | 0.999715 | 0.999853 |
| | 6 | 0.999701 | 0.999857 | 1 | 0.999812 | 0.999042 | 0.999979 | 0.999675 | 0.999182 | 0.997212 | 0.999984 | 0.999953 | 0.99715 | 0.999296 |
| | 7 | 0.999754 | 0.99994 | 1 | 0.999952 | 0.999883 | 0.999998 | 0.999998 | 0.999875 | 0.999637 | 0.999998 | 0.999998 | 0.999752 | 0.999899 |
| | 8 | 0.999579 | 0.998769 | 0.999989 | 0.998559 | 0.980343 | 0.999264 | 0.977615 | 0.996147 | 0.94226 | 0.999865 | 0.997964 | 0.939424 | 0.985815 |
| | 9 | 0.999673 | **0.999946** | **1** | 0.999936 | **0.99997** | **1** | **0.999999** | **0.999979** | **0.999953** | **1** | **1** | **0.999957** | **0.999951** |
| 3 | 0 | 0.997963 | 0.999078 | 0.982713 | 0.999367 | 0.998115 | 0.999967 | 0.990083 | 0.994797 | 0.987132 | 0.996593 | 0.998744 | 0.987148 | 0.994308 |
| | 1 | 0.998412 | 0.998405 | 0.999967 | 0.994694 | 0.989406 | 0.999908 | 0.99938 | 0.99783 | 0.992692 | 0.996846 | 0.99965 | 0.996275 | 0.996955 |
| | 2 | **0.999766** | 0.999874 | 1 | 0.999932 | 0.999854 | 0.999999 | 0.999997 | 0.999937 | 0.999747 | 0.999999 | 0.999998 | 0.99975 | 0.999904 |
| | 3 | 0.99971 | **0.999958** | **1** | **0.99996** | **0.999966** | **1** | **0.999999** | **0.999979** | **0.999946** | **1** | **1** | **0.999947** | **0.999955** |
| | 4 | 0.999738 | 0.999947 | 1 | 0.999947 | 0.999956 | 0.999999 | 0.999996 | 0.999942 | 0.999815 | 0.999999 | 0.999999 | 0.999697 | 0.99992 |
| | 5 | 0.999759 | 0.999926 | 1 | 0.999899 | 0.999909 | 0.999998 | 0.999998 | 0.999886 | 0.999595 | 0.999998 | 0.999991 | 0.999777 | 0.999895 |
| | 6 | 0.999721 | 0.999901 | 1 | 0.99983 | 0.999239 | 0.999859 | 0.999859 | 0.999391 | 0.998258 | 0.999991 | 0.999971 | 0.998011 | 0.999522 |
| | 7 | 0.999755 | 0.999956 | 1 | 0.999956 | 0.999901 | 0.999998 | 0.999998 | 0.999893 | 0.999713 | 0.999997 | 0.999997 | 0.999783 | 0.999912 |
| | 8 | 0.999613 | 0.999201 | 0.999993 | 0.998571 | 0.982176 | 0.999536 | 0.991803 | 0.99642 | 0.966812 | 0.999865 | 0.998793 | 0.954321 | 0.990592 |
| | 9 | 0.999692 | 0.999946 | **1** | 0.999946 | **0.999972** | **1** | **0.999999** | 0.999979 | **0.999954** | **1** | **1** | **0.999957** | 0.999954 |

**Table 9.** Average based on levels of each factor

| | | stage | | | |
|---|---|---|---|---|---|
| | Try | 0 | 1 | 2 | 3 |
| | 0 | 0.9997173 | 0.9997173 | 0.9997173 | 0.9997173 |
| | 1 | 0.9995158 | 0.9995158 | 0.9995158 | 0.9995158 |
| | 2 | 0.9998801 | 0.9999282 | 0.9999394 | 0.9999513 |
| | 3 | 0.9999575 | 0.999961 | 0.9999663 | 0.9999669 |
| | 4 | 0.9998977 | 0.9999439 | 0.9999515 | 0.999956 |
| | 5 | 0.9998655 | 0.9999006 | 0.9999205 | 0.9999275 |
| | 6 | 0.9991318 | 0.9996859 | 0.9997202 | 0.9998347 |
| | 7 | 0.9998328 | 0.999912 | 0.9999355 | 0.9999391 |
| | 8 | 0.9883984 | 0.9952316 | 0.9969596 | 0.9977304 |
| $c_g$ | 1 | 0.9997044 | 0.9997204 | 0.9997242 | 0.9997281 |
| | 2 | **0.9999069** | **0.9999352** | **0.9999461** | **0.9999501** |
| | 3 | 0.9957876 | 0.9982765 | 0.9988718 | 0.9991681 |
| $c_p$ | 1 | 0.9996022 | 0.9997881 | **0.9998013** | **0.9998396** |
| | 2 | **0.9997488** | **0.9997906** | 0.9998009 | 0.9998036 |
| | 3 | 0.9960480 | 0.9983534 | 0.9989398 | 0.9992031 |
| $c_w$ | 1 | **0.9998052** | 0.9998433 | 0.9998578 | 0.9998613 |
| | 2 | 0.9959572 | 0.9982361 | 0.9988139 | 0.9990710 |
| | 3 | 0.9996365 | **0.9998526** | **0.9998703** | **0.9999140** |

**Table 10.** Final decision for parameter setting.

| stage | Factor 1 | Factor 2 | Factor 3 | $c_g$ | $c_p$ | $c_w$ | $C_g$ | $C_p$ | $C_w$ |
|---|---|---|---|---|---|---|---|---|---|
| 0 | 2 | 2 | 1 | 0.4 | 0.25 | 0.3 | 0.40 | 0.65 | 0.95 |
| 1 | 2 | 2 | 3 | 0.4 | 0.25 | 0 | 0.40 | 0.65 | 0.65 |
| 2 & 3 | 2 | 1 | 3 | 0.4 | 0.35 | 0 | 0.40 | 0.75 | 0.75 |

To further verify the effectiveness of the proposed SS3OA, we performed 25 additional runs for each test problem, including a new one, that is, Try 9, based on the setting listed in Table 10. The results shown in Tables 7 and 8 confirm that Try 3 is still the best among tries 0–8. However, Try 9, based on the setting listed in Table 9, outperformed Try 3 for all stages except stage 0, even though the difference can only be determined from the eight digits after the decimal point (see Table 11).

Hence, Try 9 is adapted and compared with other state-of-the-art algorithms in Ex 2. In addition, the setting is used in the algorithms implemented in Eq. (7), when comparing with the proposed BAT-SSOA3.

**Table 11.** Average of the best $G$ for each test problem in each stage for 30 runs.

| stage | Try | 1 | 2 | 3 | 4 | 5 | 6 | 8 | 9 | 10 | 11 | 12 | 15 | average |
|---|---|---|---|---|---|---|---|---|---|---|---|---|---|---|
| 0 | 0 | 0.999553 | 0.999787 | 1 | 0.999818 | 0.999871 | 0.999995 | 0.999936 | 0.999754 | 0.999743 | 0.999999 | 0.999999 | 0.999745 | 0.9998499 |
| | 1 | 0.999599 | 0.99989 | 1 | 0.999603 | 0.999553 | 0.999996 | 0.999995 | 0.999799 | 0.999237 | 0.999996 | 0.999995 | 0.999371 | 0.9997528 |
| | 2 | 0.999781 | 0.999959 | **1** | 0.999962 | 0.999902 | 0.999999 | 1 | 0.999957 | 0.99991 | 0.999999 | 0.999999 | 0.999892 | 0.9999467 |
| | 3 | 0.999775 | **0.999967** | **1** | **0.999975** | **0.99999** | **1** | 1 | **0.999984** | 0.999973 | 1 | **1** | 0.999974 | **0.9999698** |
| | 4 | 0.99978 | 0.999963 | 1 | 0.999961 | 0.999915 | 0.999999 | 0.999999 | 0.999936 | 0.999763 | 0.999999 | 0.999999 | 0.999762 | 0.9999229 |
| | 5 | **0.999783** | 0.999945 | 1 | 0.999942 | 0.999945 | 0.999999 | 0.999998 | 0.999897 | 0.999689 | 0.999999 | 0.999999 | 0.999833 | 0.9999191 |
| | 6 | 0.999735 | 0.999866 | 1 | 0.99991 | 0.999622 | 0.999995 | 0.999945 | 0.999641 | 0.999406 | 0.999993 | 0.999983 | 0.998001 | 0.9995913 |
| | 7 | 0.99978 | 0.99995 | 1 | 0.999943 | 0.999924 | 0.999998 | 0.999998 | 0.999936 | 0.999671 | 0.999999 | 0.999998 | 0.999828 | 0.9999187 |
| | 8 | 0.999636 | 0.999703 | 0.999998 | 0.999707 | 0.995993 | 0.999937 | 0.994148 | 0.999157 | 0.997015 | 0.999979 | 0.999801 | 0.985417 | 0.9975409 |
| | 9 | 0.999779 | **0.999969** | **1** | 0.999973 | **0.999992** | **1** | 1 | **0.999987** | 1 | **1** | **1** | 0.999974 | **0.9999704** |
| 1 | 0 | 0.999553 | 0.999787 | 1 | 0.999818 | 0.999871 | 0.999995 | 0.999936 | 0.999754 | 0.999743 | 0.999999 | 0.999999 | 0.999745 | 0.9998499 |
| | 1 | 0.999599 | 0.99989 | 1 | 0.999603 | 0.999553 | 0.999996 | 0.999995 | 0.999799 | 0.999237 | 0.999996 | 0.999995 | 0.999371 | 0.9997528 |
| | 2 | 0.999783 | 0.999967 | **1** | 0.999964 | 0.999922 | 0.999999 | 1 | 0.99997 | 0.999922 | 1 | 0.999999 | 0.999999 | 0.9998984 |
| | 3 | 0.999775 | **0.999972** | **1** | **0.999975** | **0.999991** | **1** | 1 | **0.999987** | **0.999974** | **1** | **1** | **0.999976** | **0.9999711** |
| | 4 | 0.99978 | 0.99997 | 1 | 0.999967 | 0.999977 | 0.999999 | 0.999999 | 0.999969 | 0.999885 | 1 | 0.999999 | 0.999922 | 0.9999556 |
| | 5 | 0.999783 | 0.999961 | 1 | 0.999965 | 0.99997 | 0.999999 | 0.999999 | 0.999934 | 0.99988 | 1 | 0.999999 | 0.999999 | 0.9999488 |
| | 6 | 0.999767 | 0.999933 | 1 | 0.99991 | 0.999623 | 0.999995 | 0.999996 | 0.999868 | 0.99906 | 0.999996 | 0.999991 | 0.999333 | 0.9997894 |
| | 7 | **0.999783** | 0.999963 | 1 | 0.99996 | 0.999929 | 0.999999 | 0.999999 | 0.999945 | 0.999829 | 0.999999 | 0.999999 | 0.999828 | 0.9999362 |
| | 8 | 0.999709 | 0.99984 | 0.999998 | 0.999724 | 0.999284 | 0.999979 | 0.999816 | 0.999815 | 0.997015 | 0.999979 | 0.999835 | 0.997129 | 0.9992984 |
| | 9 | 0.999779 | **0.999971** | **1** | 0.999974 | **0.999993** | **1** | 1 | **0.999989** | 0.999973 | **1** | **1** | **0.999976** | **0.9999716** |
| 2 | 0 | 0.999553 | 0.999787 | 1 | 0.999818 | 0.999871 | 0.999995 | 0.999936 | 0.999754 | 0.999743 | 0.999999 | 0.999999 | 0.999745 | 0.9998498 |
| | 1 | 0.999599 | 0.99989 | 1 | 0.999603 | 0.999553 | 0.999996 | 0.999995 | 0.999799 | 0.999237 | 0.999996 | 0.999995 | 0.999371 | 0.9997527 |
| | 2 | 0.999783 | 0.999967 | **1** | 0.999967 | 0.99998 | 1 | 1 | 0.999973 | 0.999942 | 1 | 0.999999 | 0.999918 | 0.9999604 |
| | 3 | 0.999775 | **0.999972** | **1** | **0.999975** | **0.999993** | **1** | 1 | **0.999987** | **0.99998** | **1** | **1** | **0.999976** | **0.9999719** |
| | 4 | 0.99978 | 0.999972 | 1 | 0.999968 | 0.999982 | 1 | 0.999977 | 0.999911 | 1 | 0.999999 | 0.999999 | 0.999922 | 0.9999591 |
| | 5 | 0.999783 | 0.999961 | 1 | 0.99997 | 0.99997 | 0.999999 | 0.999999 | 0.999934 | 0.999893 | 1 | 0.999999 | 0.99988 | 0.9999488 |
| | 6 | 0.999775 | 0.999933 | 1 | 0.99991 | 0.999739 | 0.999997 | 0.999996 | 0.999868 | 0.999311 | 0.999998 | 0.999994 | 0.999333 | 0.9998210 |
| | 7 | **0.999784** | 0.999963 | 1 | 0.999964 | 0.999955 | 0.999999 | 0.999999 | 0.999945 | 0.999854 | 0.999999 | 0.999999 | 0.999883 | 0.9999455 |
| | 8 | 0.999753 | 0.999814 | 0.999999 | 0.999726 | 0.999284 | 0.999979 | 0.999816 | 0.999375 | 0.997015 | 0.999985 | 0.999876 | 0.997129 | 0.9993126 |
| | 9 | 0.999779 | **0.999976** | **1** | 0.999974 | **0.999993** | **1** | 1 | **0.999989** | 0.999974 | **1** | **1** | **0.999979** | **0.9999723** |
| 3 | 0 | 0.999553 | 0.999787 | 1 | 0.999818 | 0.999871 | 0.999995 | 0.999936 | 0.999754 | 0.999743 | 0.999999 | 0.999999 | 0.999745 | 0.9998499 |
| | 1 | 0.999599 | 0.99989 | 1 | 0.999603 | 0.999553 | 0.999996 | 0.999995 | 0.999799 | 0.999237 | 0.999996 | 0.999995 | 0.999371 | 0.9997528 |
| | 2 | 0.999783 | 0.999968 | **1** | 0.999968 | 0.99998 | 1 | 1 | 0.999975 | 0.999942 | 1 | 1 | 0.999918 | 0.9999608 |
| | 3 | 0.999778 | 0.999972 | **1** | **0.999975** | **0.999994** | **1** | 1 | **0.999987** | **0.999981** | **1** | **1** | **0.999977** | **0.9999724** |

| | | | | | | | | | | | | | | | | | |
|---|---|---|---|---|---|---|---|---|---|---|---|---|---|---|---|---|---|
| 4 | 0.999782 | **0.999972** | **1** | 0.999969 | 0.999983 | | **1** | | **1** | 0.999981 | 0.999927 | | **1** | | **1** | 0.999922 | 0.9999614 |
| 5 | **0.999784** | 0.999964 | 1 | 0.999967 | 0.99997 | 0.999999 | 1 | | 0.999941 | 0.999893 | | 1 | 0.999999 | | 0.99992 | 0.9999531 |
| 6 | 0.999781 | 0.999933 | | 1 | 0.999912 | 0.999997 | 0.999999 | 0.999868 | 0.999714 | 0.999998 | 0.999994 | 0.999711 | 0.9998956 |
| 7 | 0.999784 | 0.999967 | | 1 | 0.999971 | 0.999972 | 0.999999 | | 1 | 0.999963 | 0.999888 | | 1 | | 1 | 0.999927 | 0.9999559 |
| 8 | 0.999776 | 0.999814 | 0.999999 | 0.999886 | 0.999284 | 0.999979 | 0.999816 | 0.999375 | 0.997015 | 0.999996 | 0.999876 | 0.997867 | 0.9993903 |
| 9 | **0.999784** | **0.999976** | | **1** | **0.999976** | 0.999993 | | **1** | | **1** | **0.999989** | 0.999976 | | **1** | | **1** | **0.999979** | **0.9999731** |

## 4.3 Ex 2

There are no existing algorithms for the proposed general RRAP. To demonstrate its performance in Ex 2, we modified well-known and important traditional RRAP-related algorithms by adding connected vectors obtained from the BAT to calculate the reliability, that is, BAT-SSO, BAT-BSO, BAT-nSSO, , BAT-ifSSO, BAT-PSO, and BAT-PSSO.

The update mechanism of each algorithm used for comparison in Ex 2 is as follows:

UM0: Eq. (1).

UM1:

$$x_j^* = \begin{cases} g_j & \text{if } \rho_{[0,1]} \in [0., C_g) \\ p_j & \text{if } \rho_{[0,1]} \in [C_g, C_p) \\ x_j & \text{if } \rho_{[0,1]} \in [C_p, C_w) \\ x & \text{if } \rho_{[0,1]} \in [C_w, 1.] \end{cases} + \begin{cases} 0 & \text{if } x_j \text{ is a discrete variable or } \rho_{[0,1]} \geq C_w \\ l & \text{otherwise} \end{cases}, \tag{12}$$

$$l = 0.0005 \cdot \rho_{[-0.5,0.5]} \cdot \frac{\text{gen}}{genBest}. \tag{13}$$

UM2: The boundary update as discussed in Section 3.3.

UM3: Similar to UM1, but replacing Eq. (13) with Eq. (14).

$$l = \rho_{[0,1]} \cdot \exp(\frac{-100 \cdot \text{gen}}{N_{\text{gen}}}). \tag{14}$$

UM4:

1) Randomly select and reinitialize two floating-point variables in $G$.

2) Randomly select one floating-point variable, for example $r_j$, and reset the value of $r_j$ based on Eq. (5).

UM5 (for PSO):

$$V_i = \begin{cases} V_{\text{ub}} & \text{if } [wV_i + c_1 r_{[0,1]}(P_i - X_i) + c_2 r_{[0,1]}(G - X_i)] > V_{\text{ub}} \\ V_{\text{lb}} & \text{if } [wV_i + c_1 r_{[0,1]}(P_i - X_i) + c_2 r_{[0,1]}(G - X_i)] < V_{\text{lb}} \\ wV_i + c_1 r_{[0,1]}(P_i - X_i) + c_2 r_{[0,1]}(G - X_i) & \text{otherwise} \end{cases} \quad (15)$$

where $w = 0.9$, $c_1 = c_2 = 2$.

$$X_i = \begin{cases} X_{\text{ub}} & \text{if } (X_i + V_i) > X_{\text{ub}} \\ X_{\text{lb}} & \text{if } (X_i + V_i) < X_{\text{lb}} \\ X_i + V_i & \text{otherwise} \end{cases} \quad (16)$$

Summary of the updated mechanism for each algorithm is listed in Table 12.

**Table 12.** Updated mechanism of each algorithm.

| Algorithms | gBest | Not gBest | Integer variables | Floating variables | All variables |
|---|---|---|---|---|---|
| BAT-SSOA3 | UM2 | UM0 | | | |
| BAT-SSO | | | | | UM0 |
| BAT-BSO | UM4 | UM0 | | | |
| BAT-nSSO | | | UM0 | UM1 | |
| BAT-ifSSO | | | UM0 | UM3 | |
| BAT-PSSO | | | UM0 | UM5 | |
| BAT-PSO | | | | | UM5 |

**Table 13.** Average of the best $G$ for each test problem, stage, and algorithm for 30 runs.

| S* | A* | P*1 | P*2 | P*3 | P*4 | P*5 | P*6 | P*8 | P*9 | P*10 | P*11 | P*12 | P*13 | average |
|---|---|---|---|---|---|---|---|---|---|---|---|---|---|---|
| 0 | s | 0.999767 | **0.999972** | 1 | 0.999971 | 0.999991 | **1** | **1** | 0.999985 | 0.999969 | **1** | 1 | 0.999967 | 0.9999685 |
| | g | 0.999768 | 0.999827 | 1 | 0.999904 | 0.999862 | 1 | 1 | 0.999958 | 0.999744 | 1 | 0.999999 | 0.999796 | 0.99990483 |
| | i | 0.999626 | 0.999872 | 1 | 0.999816 | 0.999787 | 0.999997 | 0.999997 | 0.999842 | 0.999316 | 0.999993 | 0.999987 | 0.999758 | 0.99983258 |
| | n | 0.999776 | 0.999963 | 1 | 0.999956 | 0.999958 | 0.999999 | 0.999999 | 0.999945 | 0.999827 | 1 | 0.999999 | 0.999849 | 0.99993925 |
| | p | 0.999316 | 0.992367 | 0.999913 | 0.969406 | 0.999188 | 0.99957 | 0.999188 | 0.984995 | 0.990263 | 0.999396 | 0.999342 | 0.978414 | 0.99253175 |
| | o | 0.999723 | 0.999911 | 1 | 0.99994 | 0.999908 | 1 | 1 | 0.999903 | 0.999774 | 0.999998 | 0.999999 | 0.99988 | 0.99991967 |
| | a | **0.999779** | 0.999969 | 1 | **0.999973** | **0.999992** | 1 | 1 | **0.999987** | **0.99997** | 1 | 1 | **0.999974** | **0.99997042** |
| 1 | s | 0.999777 | **0.999972** | 1 | 0.999972 | 0.999992 | 1 | 1 | 0.999985 | **0.999976** | 1 | 1 | **0.99998** | 0.99997117 |
| | g | 0.999759 | 0.999945 | 1 | 0.999906 | 0.999937 | 1 | 1 | 0.999934 | 0.999867 | 0.999999 | 1 | 0.999898 | 0.99993708 |
| | i | 0.999683 | 0.999872 | 1 | 0.999816 | 0.999787 | 0.999998 | 0.999999 | 0.999895 | 0.999316 | 0.999995 | 0.999998 | 0.999758 | 0.99984308 |
| | n | 0.999776 | 0.999967 | 1 | 0.999964 | 0.999988 | 1 | 1 | 0.999975 | 0.999949 | 1 | 1 | 0.999996 | 0.99996475 |
| | p | 0.999316 | 0.994956 | 0.999913 | 0.992612 | 0.972539 | 0.998027 | 0.999337 | 0.970658 | 0.958064 | 0.996443 | 0.999334 | 0.980504 | 0.98847525 |
| | o | 0.999763 | 0.999941 | 1 | 0.999948 | 0.999967 | **1** | 1 | 0.999939 | 0.999836 | 1 | 0.999999 | 0.999892 | 0.9999405 |
| | a | **0.999774** | 0.999971 | 1 | **0.999974** | **0.999993** | 1 | 1 | **0.999989** | 0.999973 | 1 | 1 | **0.999979** | **0.99997158** |
| 2 | s | **0.99998** | 0.999972 | 1 | 0.999973 | 0.999992 | 1 | 1 | 0.999985 | **0.999978** | 1 | 1 | **0.999981** | 0.99997175 |
| | g | 0.999776 | 0.999959 | 1 | 0.999856 | 0.999936 | 1 | 1 | 0.999881 | 0.999727 | 0.999998 | 1 | 0.999908 | 0.99992008 |
| | i | 0.999688 | 0.999922 | 1 | 0.999898 | 0.99998 | 0.999998 | 0.999999 | 0.999903 | 0.999316 | 0.999998 | 0.999998 | 0.999758 | 0.99985608 |
| | n | **0.99998** | 0.999967 | 1 | 0.999967 | 0.999989 | 1 | 1 | 0.999983 | 0.999954 | 1 | 1 | 0.99996 | 0.99996667 |
| | p | 0.999316 | 0.994956 | 0.999913 | 0.992612 | 0.972539 | 0.998027 | 0.999121 | 0.995891 | 0.978625 | 0.996443 | 0.987458 | 0.980504 | 0.99128375 |
| | o | 0.999769 | 0.999941 | 1 | 0.999952 | 0.999967 | 1 | 1 | 0.999961 | 0.999934 | 1 | 0.999999 | 0.999892 | 0.9999515 |
| | a | 0.999779 | **0.999976** | 1 | **0.999974** | **0.999993** | 1 | 1 | **0.999989** | 0.999974 | 1 | 1 | 0.999979 | **0.99997233** |
| 3 | s | 0.999781 | 0.999973 | 1 | **0.999976** | 0.999992 | 1 | 1 | 0.999987 | **0.999979** | 1 | 1 | **0.999981** | 0.99997242 |
| | g | 0.999756 | 0.999742 | 1 | 0.999825 | 0.999863 | 1 | 1 | 0.999943 | 0.999832 | 0.999993 | 0.999999 | 0.999898 | 0.99990425 |
| | i | 0.999708 | 0.999922 | 1 | 0.999898 | 0.9998 | 0.999998 | 0.999999 | 0.999903 | 0.999467 | 0.999999 | 0.999998 | 0.999796 | 0.999874 |
| | n | 0.999783 | 0.999969 | 1 | 0.999969 | 0.999989 | 1 | **1** | 0.999985 | 0.999957 | **1** | 1 | 0.99997 | 0.99996867 |
| | p | 0.999316 | 0.994956 | 0.999957 | 0.992612 | 0.980807 | 0.998027 | 0.999121 | 0.975144 | 0.978625 | 0.999774 | 0.987458 | 0.980504 | 0.99052508 |
| | o | 0.999769 | 0.999941 | 1 | 0.999953 | 0.999976 | 1 | 1 | 0.999961 | 0.999934 | 1 | 0.999999 | 0.999902 | 0.99995317 |
| | a | **0.999784** | **0.999976** | 1 | **0.999976** | **0.999993** | 1 | 1 | **0.999989** | 0.999974 | **1** | 1 | 0.999979 | **0.99997308** |

S\*: stage
A\*: algorithm
P\*: problem

The result of Ex 2 contrasts the best final *gBest*s obtained from the proposed BAT-SSOA3 with the best parameter setting obtained from Ex 1. In Ex 2, we focus on the comparisons among the best solutions obtained from all the abovementioned algorithms.

BAT-SSO, BAT-gSSO, BAT-ifSSO, BAT-nSSO, BAT-PSO, BAT-PSSO, and BAT-SSOA3 are assigned as algorithms *s*, *g*, *i*, *n*, *p*, *o*, and *a*, respectively. The obtained average of the best gBests is shown in Table 13, and the best among all the algorithms is written in bold. The first row of Table 13 lists stages 0–3; the second row lists the algorithms, that is, *s*, *g*, *i*, *n*, *p*, *o*, and *a*; columns 3–15 indicate the averages of the best gBests for problems 1–12, respectively; and the last column lists the average of the values from columns 3 to 15.

To evaluate the performance, including the efficiency and effectiveness (i.e., solution quality) of the proposed BAT-SSOA3, as in Ex 1, and these performed in all existing algorithms for the RRAP, we only provided the best of the final *gBest*s for each problem and algorithm in Table 13 [1, 2, 3, 4, 5, 6].

Based on the average of the best final gBest averages, that is, the last column in Table 13, the proposed algorithm BAT-SSOA3, that is, algorithm *a*, showed better performance than the others for all stages in all cases. In addition, the proposed algorithm outperformed the other algorithms for all stages except problems 1, 10, and 15.

Surprisingly, BAT-SSO, which simply used Eq. (7) to update variables after using the affixation solution structure, was the second best. Hence, Eq. (7) is also powerful for floating-point variables; however, this is not consistent with the conclusions found in the literature. The main reason is that the SSO parameters, that is, $C_g$, $C_p$, and $C_w$, are all adapted from Try 9, which is the best among Try 0–9 in the proposed SS3OA, as discussed in Ex 1. In addition, the results obtained from Try 3 in Ex 1 are better than those of SSO. Therefore, this observation further confirms that the proposed SS3OA in BAT-SSOA3 is powerful for tuning the parameters and selected items in Eq. (1).

Based on the data listed in Table 13, we have the following observations:

1.  In general, BAT-SSOA3 >> BAT-SSO >> BAT-nSSO >> BAT-gSSO >> BAT-PSSO >> BAT-

ifSSO >> BAT-PSO, where A >> B means that A is better than B in terms of solution quality.

2. The difference between BAT-SSOA3 and BAT-SSO is the boundary update. Hence, the boundary is helpful in improving the solution quality.

3. The differences among BAT-SSO, BAT-nSSO, BAT-ifSSO, and BAT-PSSO is the way the floating-point variables are updated. As BAT-SSO >> BAT-nSSO >> BAT-PSSO >> BAT-ifSSO, UM 0 >> UM 1 >> UM 5 >> UM 3 in the updated floating-point variables.

4. The difference between BAT-PSSO and BAT-PSO is the updated integer variable. Hence, UM 0 >> UM 5 in the updated integer variable as BAT-PSSO >> BAT-PSO.

5. As BAT-SSOA3 >> BAT-BSO, and the difference is the way gBest is updated, UM 2 >> UM 4 in the updated gBest.

Hence, based on the discussion above for both Ex 1 and Ex 2, the proposed BAT-SSOA3 with parameter settings from Ex 1 can balance both global and local searches for improving solution quality compared with previously published RRAP algorithms.

## 5. Conclusions

The serial-parallel or bridge structure RRAP is not practical and reasonable for all types of applications in real-life networks. Thus, a novel RRAP called GRRAP is proposed to address this problem. A new BAT-SSOA3, integrated BAT-connected vectors, SS3OA, boundary update, penalty function, and affixation solution structure are used to solve the proposed GRRAP in this study.

The proposed BAT-SSOA3 is highly efficient and effective. One of the major difficulties in the proposed general RRAP is the calculation of its objective function, that is, the reliability of a general binary-state network. However, after the implementation of the BAT algorithm, the difficult binary-state network reliability can be easily solved by adding the probabilities of all connected vectors without tedious and complicated calculation on reliability, which is the objective function in the GRRAP. Moreover, the proposed BAT-connected vector can be implemented in any existing algorithm as an analytical tool or function to calculate all types of binary-state networks.

The SS3OA can be systematically and efficiently used to obtain the best parameters and the most suitable updated mechanism for Eq. (7). The use of the boundary condition also provides a significant improvement in the solution quality after comparison with other algorithms. In addition, the affixation solution structure combines a pair of different types of variables into one variable to reduce the runtime by half. The penalty function avoids boundary obstacles occurring in the feasible solution space and strengthens the ability of SSO variants to solve problems.

Based on the fair and complete comparisons of experimental results from Ex 2, as discussed in Section 4.2, the solution quality of BAT-SSOA3 outperforms and is superior to well-known related methods, which use BAT to calculate reliability. The penalty function avoids boundary obstacles and the affixation solution structure reduces the update time by half. In future studies, the proposed BAT-SSOA3 will be extended to warm-standby RRAP [38], cold-standby RRAP [2, 39-41], multi-objective RRAP [5, 42], multi-state RRAP [43-46], and different optimization problems with more variables or larger-scale benchmarks.

## Acknowledgements

This research was supported in part by the Ministry of Science and Technology, R.O.C. under grant MOST 107-2221-E-007-072-MY3 and MOST 104-2221-E-007-061-MY3.